\documentclass[12pt,draftcls,onecolumn]{IEEEtran}
\usepackage{amssymb,latexsym}
\usepackage{graphicx}

\newcommand{\argmax}{\mathop{\mbox{\rm arg\,max}}}
\newcommand{\argmin}{\mathop{\mbox{\rm arg\,min}}}
\newtheorem{thm}{Theorem}[section]

\begin{document}
\sloppy

\title{An Index-based Deterministic Asymptotically Optimal Algorithm for Constrained Multi-armed Bandit Problems}
\author{Hyeong Soo Chang
\thanks{H.S. Chang is with the Department of Computer Science and Engineering at Sogang University, Seoul 121-742, Korea. (e-mail:hschang@sogang.ac.kr).}%
}

\maketitle
\begin{abstract}
For the model of constrained multi-armed bandit,
we show that by construction there exists an index-based deterministic 
asymptotically optimal algorithm. 
The optimality is achieved by the convergence of the probability of 
choosing an optimal feasible arm to one over infinite horizon.
The algorithm is built upon Locatelli \emph{et al.}'s ``anytime 
parameter-free thresholding" algorithm under the assumption that
the optimal value is known. We provide a finite-time bound to the 
probability of the asymptotic optimality given as $1-O(|A|Te^{-T})$ 
where $T$ is the horizon size and $A$ is the set of the arms in the bandit.
We then study a relaxed-version of the algorithm in a general form that estimates 
the optimal value and discuss the asymptotic optimality of the algorithm
after a sufficiently large $T$ with examples.
\end{abstract}


\section{Introduction}

Consider a constrained multi-armed bandit (CMAB)~\cite{changcmab} problem where there is a finite set $A$ of arms, $|A|\geq 2$, and a single arm in $A$ needs to be sequentially played. 
When $a$ in $A$ is played at discrete time $t\geq 1$,
the player not only obtains a sample bounded reward $X_{a,t}\in \Re$ 
drawn from an unknown reward-distribution associated with $a$, whose unknown finite expectation and finite variance are $\mu_a$ and $\sigma_{R,a}^2$, respectively, 
but also obtains a sample bounded cost $Y_{a,t} \in \Re$ drawn from an unknown cost-distribution associated with $a$, whose unknown expectation and variance are $C_a$ and $\sigma_{C,a}^2$, respectively. 
Sample rewards and costs across arms are all independent for all time steps. That is, 
$X_{a,t}, X_{b,s}, Y_{p,t'}$, and $Y_{q,s'}$ are independent for all $a,b,p,q\in A$ and all $t,s,t',s'\geq 1$.
For any fixed $a$ in $A$, $X_{a,t}$'s and $Y_{a,t}$'s for $t\geq 1$ are identically distributed, respectively.
We define the \emph{feasible set} $A_f$ of arms such that $A_f := \{ a\in A | C_a \leq C \}$ for a constant $C \in \Re$ known to the player and assume that $A_f \neq \emptyset$.
Our goal is to find an \emph{optimal feasible} arm that achieves \emph{the optimal value} $\mu^*:=\max_{a\in A_f}\mu_a$.
(For the sake of simplicity, we consider one constraint case. It is straightforward to extend our results into multiple-constraints case.)

The model of unconstrained MAB has been used for studying many (practical) problems
(see, e.g.,~\cite{berry}~\cite{gittins}~\cite{cesa} for in depth cover of the topic and examples). 
However, there also exist related MAB problems that are involved with one or more of \emph{conflicting} objective functions with the main objective functions. 
These conflicting objective functions play the roles of \emph{constraints} for optimizing the main objective functions in CMAB problems.
For example, a trade-off exists between achieving a ``small" delay (or ``high" throughput) and ``low" power consumption in wireless communication networks.
To maximize the throughput (or to minimize the delay) we need to transmit with the highest
available power level because it will increase the probability of successful transmission.
On the other hand, to minimize the power consumption, we need to transmit with the lowest power level available.
We can consider the problem of selecting an optimal feasible power level among all available powers that keeps the delay cost below some given bound. 
In fact, in many scheduling and queueing control problems, there exist certain
trade-offs between ``throughput" and ``delay" in general.
%Another quick example problem can be about choosing an optimal feasible Ad that maximizes some revenue that keeps the marketing cost below some bound in Ad placement management for some internet sites,

We define an \emph{algorithm} $\pi:=\{\pi_t, t=1,2,...\}$ as a sequence of mappings such that $\pi_t$ maps from the set of past plays and rewards and costs, $H_{t-1}:=(A \times \Re \times \Re)^{t-1}$ if $t\geq 2$ and $\emptyset$ if $t=1$, to the set $A$.
We denote the set of all possible such deterministic algorithms as $\Pi$.
The \emph{asymptotic optimality} introduced by Robbins~\cite{robbins} for the optimality in terms of \emph{transient} behavior of an algorithm
will be used as the measure of the performance.
Let $A^*_f := \{a\in A_f| \mu_a = \mu^*\}$ and let $I^{\pi}_t$ denote the arm selected by $\pi$ at time $t$.
Given $\pi$ in $\Pi$, we say that $\pi$ is asymptotically optimal
if $\sum_{a\in A^*_f} \Pr\{I^{\pi}_t = a\} \rightarrow 1$ as $t \rightarrow \infty$.

This note begins with presenting an algorithm, called ``Constrained Anytime Parameter-free Thresholding (CAPT)," in order to show that by construction there exists 
an asymptotically optimal \emph{index-based deterministic} algorithm in $\Pi$.
It has been conjectured that devising such an algorithm is difficult~\cite{changcmab} 
because the question of how to ``mix" a process of estimating the feasibility of each arm 
into the exploration (and possibly exploitation) process of estimating the reward-optimality of each feasible arm by \emph{one deterministic index} needs to be answered.
This note provides an affirmative report to the question.

To approach a CMAB problem, instead of searching a proper index for each action, one can consider a methodology that ``separates" the two associated problems of CMAB ``in time" 
by solving the cost-feasibility problem first and then solves the reward-optimality problem conditioned on the results about the feasibility, eventually achieving the asymptotic optimality. 
The immediate open questions are firstly if the method is in $\Pi$ (even with putting aside an index-based selection) 
because a strong negative argument is that
the method would provide only \emph{a probabilistic judgement} in order to change the tune and secondly if the algorithm is analyzable in terms of the asymptotic optimality.
A \emph{randomized} strategy recently studied by Chang~\cite{changcmab} 
extends the $\epsilon$-greedy MAB strategy~\cite{auer} and
works with
the underpinning exploration method of \emph{uniform random-selection}.
It is worthwhile to note that the arguments in the analysis regarding the 
asymptotic optimality of the strategy were provided  
with the very above idea of separating the estimation process in the bandit 
in time by the probabilistic estimation result for the feasibility.
This is from the ground-level process of 
the uniform selection that allows conditioning a ``guaranteed level" of
the feasibility estimation for a given time.
With controlling the values of $\epsilon$ used over time by the strategy, it was
proved that the strategy achieves the asymptotic optimality after a sufficiently 
large horizon.
The strategy is simple even if randomized but \emph{not index-based}.
Still, the strategy provides an important direction towards designing 
solution methods for CMAB problems.

In addition, a potentially profitable functionality missed by the constrained 
$\epsilon$-greedy strategy is the usage of problem-characteristics about the
reward-optimality and the cost-feasibility.
Let us define $\Delta_i^\epsilon:=|\mu_i-\mu^*|+\epsilon, i\in A$ and $\Phi_i^\epsilon:=|C_i-C|+\epsilon,i\in A$ where $\epsilon \geq 0$. 
The role of $\epsilon$ is to provide some tolerance in what we measure.
These values can guide us for determining degrees of allocating samples over arms.
The more closely competitive and feasible arms there are, the more difficult problem is in general.
For example, if the second best arm becomes very competitive with the best, the harder distinguishing the best from the second best is.
If $\Phi_i^\epsilon$ becomes closer to zero, checking the feasibility of $i$
becomes more difficult.
More sampling efforts need to be put in distinguishing closely feasible and competitive arms.
Indeed, the probabilities related with the convergences have been expressed in terms of a function of $\Delta_i^\epsilon$'s, called ``\emph{complexity} of the problem" (see, e.g.,~\cite{audi}~\cite{bubeck}~\cite{kaufmann}~\cite{locatelli}, etc.~and also confer with, e.g.,~\cite{russo} about related but different complexity measures).
Noticeably, 
Locatelli \emph{et al.}~\cite{locatelli} developed 
an index-based deterministic algorithm, called ``Anytime Parameter-free Thresholding (APT)" 
for ``thresholding bandit" problems by using $\Phi_i^\epsilon$'s.
It turns out that the problem considered exactly coincides with the cost-feasibility problem in CMAB. 
In particular, the index of $a$ in $A$ is given by $\bar{\Phi}_a^\epsilon(t) \sqrt{T_a(t)}$
where $\bar{\Phi}_a^\epsilon(t)$ is an estimate of $\Phi_a^\epsilon$ 
obtained by replacing $C_a$ by the sample mean up to time $t$ by $T_a(t)$ cost-samples of playing $a$. 
The number of times $a$ has been played up to time $t$ is denoted by $T_a(t)$.
The index measures ``to what degree $a$ needs to be sampled" at time $t\geq 1$ of 
APT and 
APT plays an arm in the argument set that achieves the minimum index.
We can see that the index has a structure that the arm selection is affected by the values of $\bar{\Phi}_a^\epsilon$ and $\sqrt{T_a(t)}$ together. 

The CAPT algorithm presented in Section~\ref{sec:algo} employs the same form of the index 
of APT but with some extension: the index of $a$ is given as
\[
K_a(t) := \min\left( \bar{\Delta}_a^\epsilon(t), \bar{\Phi}_a^\epsilon(t) \right) \sqrt{T_a(t)}.
\] The term $\bar{\Delta}_a^\epsilon(t)$ is an estimate of $\Delta_a^\epsilon$. 
Similar to $\Phi_i^\epsilon$, the sample mean of $a$ up to time $t$ by $T_a(t)$ reward-samples is used in place of $\mu_i$ in $\Delta_i^\epsilon$. 
The idea is simple.
\emph{We pose the reward-optimality problem as another cost-feasibility problem}. Each index for the two feasibility problems are combined into a new index by the \emph{minimum operator}. 
The index then measures
not only to what degree $a$ needs to be sampled for cost-feasibility but also for reward-optimality.
We prove that CAPT constructed from this simple fusion achieves the asymptotic optimality. We provide a finite-time lower bound to the probability of finding an optimal feasible action with $1-O(|A|Te^{-T})$ for a given finite horizon $T$.

CAPT works with the crucial assumption that the optimal value $\mu^*$ 
is known because $\bar{\Delta}_a^\epsilon(t),a\in A,$ needs to be computed.
However, we argue that this \emph{theoretical} study is an important step towards understanding the solvability and the complexity of CMAB. 
In fact, the procedures of some algorithms for MAB in the literature 
were given with the optimal value (or a functional value of it or a known bound to it) as an input parameter 
and accordingly analyzed
(see, e.g., Theorem 3 for the $\epsilon$-greedy algorithm in~\cite{auer}, Theorem 3 for APT in~\cite{locatelli}, Theorem 1 for UCB-E in~\cite{audi}, Theorem 2 and the related work section in~\cite{yang}, etc.).

In Section~\ref{sec:relax},
we study an algorithm in a general form, called CAPT-E (CAPT with Estimation),
where the index of CAPT is replaced by $\kappa_a(t)$ such that
\[
\kappa_a(t) := \min\left( \bar{\Delta}_a^{\epsilon,*}(t), \bar{\Phi}_a^\epsilon(t) \right) \sqrt{T_a(t)},
\] where $\bar{\Delta}_a^{\epsilon,*}(t)$ is an estimate of $\bar{\Delta}_i^\epsilon$
that substitutes $\mu^*(t)$ into $\mu^*$.
We discuss a sufficient condition that makes CAPT-E achieve the asymptotic optimality 
after a sufficiently large $T$ and some examples of $\mu^*(t)$.

The main goal of this note is to establish the \emph{existence} of an asymptotically optimal index-based deterministic algorithm in $\Pi$ and to provide a \emph{theoretical characterization} about the CMAB problems (as in Theorem 2 of Lai and Robbins~\cite{lai}).
The performance of CAPT would be a baseline for comparison or improvement
for an index-based algorithm in $\Pi$ with the criterion of the asymptotic optimality.
Finally, we show that the critical assumption can be relaxed and open some direction for further research in developing algorithms for solving CMAB problems.

\section{Related Works}

The model of CMAB is a special case of constrained Markov decision process (CMDP)~\cite{altman}~\cite{denardo}, in which we assume that all of the distributions of rewards and costs associated with all arms are \emph{unknown} to the decision maker.
Because of the assumption, the exact solution method, e.g., linear programming, is not applicable for solving CMAB problems.

Much attention has been paid recently to the model, ``Budgeted MAB (BMAB)," 
that adds a certain constraint for optimality (see, e.g.,~\cite{ding}).
In our terms, consider a random variable that takes the value of the sum of the random costs obtained by running an algorithm $\pi$ in $\Pi$ over $T$ horizon, i.e., $\sum_{t=1}^{T} Y_{I^{\pi}_t,t}$.
Let the stopping time $Q^{\pi}(B) = \min \{ T |  \sum_{t=1}^{T} Y_{I^{\pi}_t,t} > B \}$ where $B>0$ is a problem parameter called \emph{budget}.
The player stops playing at $Q^{\pi}(B)$ once it consumes up all of the budget given by $B$. 
Take the expected value of the sum of the random rewards obtained by following $\pi$ over the sample path of length $Q^{\pi}(B)-1$. We wish to maximize the expected value over all possible $\pi$. 
In other words, the goal is to obtain $\max_{\pi\in \Pi} E[\sum_{t=1}^{ Q^{\pi}(B) - 1} X_{I^{\pi}_t,t}]$ or an algorithm that achieves it. 
The key difference from CMAB is that in BMAB, the budget constraint is put on the \emph{played arm sequence}. 
Furthermore, while CMAB is a special case of CMDP as we mentioned before, it seems that BMAB is not directly related with CMDP.

Constrained simulation optimization, under the topic of ``constrained ranking and selection," considers a similar simulation setting where the values of objective and constraint functions can be obtained only by a sequential sampling process.
However, we do not draw multiple samples of reward and cost at a single time step. 
No particular assumptions on the reward and the cost distributions (e.g., normality) are made.
Sampling plan or sampling allocation is not computed in advance
as these or subset of these are common assumptions and approaches in the literature (see, e.g.,~\cite{pasupathy}~\cite{hunter}~\cite{park} and the references therein).

Various measures of studying the behaviours of the MAB algorithms exist (see, e.g., a discussion in~\cite{russo}). 
The most notable ones are probably the \emph{expected} regret~\cite{lai} for \emph{average} behaviour and the asymptotic optimality for transient behaviour.
Auer \emph{et al.}~\cite{auer} relates the asymptotic optimality with ``instantaneous" regret given as $\sum_{a\in A\setminus A^*_f} \Pr\{I^{\pi}_t = a\}$ and note that the instantaneous regret is a stronger notion than the expected regret in the convergence. 
The asymptotic optimality is directly related with
the probability of \emph{identifying} a best arm~\cite{audi}~\cite{bubeck}~\cite{kaufmann} in the so-called ``pure exploration" problem.
In the simulation optimization literature, the probability has been often referred to as
the probability of correct selection (see, e.g.,~\cite{fu}, etc.).
The different reference for the probability seems to depend on the context of the problem topic under study in the relevant literature.

The literature in MAB has rather focused on the expected regret since the work of Lai and Robbins~\cite{lai} and more particularly since Auer \emph{et al.}'s finite-time analysis on index-based algorithms~\cite{auer} (see, e.g.,~\cite{cesa} and the references therein).
It is difficult to find a work that studies the instantaneous behaviour of the existing MAB algorithms designed for the expected regret, e.g., UCB~\cite{auer} or its variants~\cite{cesa}.
That is, even if the expected behaviour of an algorithm \emph{relative to the best algorithm} has been extensively studied in the literature, the expected behavior of the algorithm itself seems to be not known yet.
Note that obtaining the expected behavior of $\sum_{a\in A} \Pr\{I^{\mbox{UCB}}_{t}=a\}, t<\infty$, for UCB essentially requires analyzing the transient behavior of UCB, i.e., the probability of $\Pr\{I^{\mbox{UCB}}_{t}=a\}$. 

Defining the expected regret within CMAB is not straightforward.
If we try a definition given by the expected loss relative to the cumulative expected reward of taking an optimal feasible arm due to the fact that the algorithm does not always play an optimal feasible arm, 
$\mu^*T-\sum_{a\in A} \mu_a (\sum_{t=1}^T \Pr\{I^{\pi}_t =a \})$ for $T$ in $[1,\infty)$, the loss can be negative.
The problem of minimizing the regret is no longer meaningful because 
this is like having a negative cycle in a shortest-path problem.
In some cases, the minimum is simply achieved by an algorithm that always plays
an infeasible arm whose reward average is higher than $\mu^*$. A possible
leverage would be introducing a function over $A$ that penalizes
an infeasibility to some degree inside the summation. 
Defining the expected ``regret" 
and design and analysis of proper algorithms will depend on the definition.
The study on the expected regret in CMAB is beyond the scope of this note and 
is left as a future research.

\section{Constrained APT Algorithm}
\label{sec:algo}
\subsection{Algorithm}

Once $I^{\pi}_t$ in $A$ is played by CAPT (referred to as $\pi$ wherever possible) at time $t$, 
a sample reward of $X_{I^{\pi}_t,t}$ and a sample cost of $Y_{I^{\pi}_t,t}$ are obtained independently.
We let $T_a(t) := \sum_{n=1}^t [ I^{\pi}_n = a ]$ where $[\cdot]$ denotes the indicator function, i.e., $[ I^{\pi}_n = a ]=1$ if $I^{\pi}_n=a$ and 0 otherwise.
The sample average-reward $\bar{X}_{T_a(t)}$ for $a$ in $A$ is then given such that $\bar{X}_{T_a(t)} = \frac{1}{T_a(t)} \sum_{n=1}^{t} X_{a,n} [ I^{\pi}_n = a ]$ 
if $T_a(t)\geq 1$ and 0 otherwise,
Similarly, $\bar{Y}_{T_a(t)}$ for $a$ in $A$ is given such that $\bar{Y}_{T_a(t)} = \frac{1}{T_a(t)} \sum_{n=1}^{t} Y_{a,n} [ I^{\pi}_n = a ]$ if $T_a(t)\geq 1$ and 0 otherwise. 
Note that $E[X_{a,t}]=\mu_a$ and $E[Y_{a,t}]=C_a$ for all $t$.
We let $\bar{\Delta}_i^\epsilon(t) = |\bar{X}_{T_i(t)} - \mu^*| + \epsilon$
and $\bar{\Phi}_i^\epsilon(t) = |\bar{Y}_{T_i(t)} - C| + \epsilon$ for $\epsilon \geq 0.$
A pseudocode for CAPT is provided below.
\\
\noindent\textbf{The Constrained APT (CAPT) algorithm $\pi$}
\begin{itemize}
\item[1.] \textbf{Initialization:} 
\begin{itemize}
\item[1.1] Select $\epsilon \geq 0$.
\item[1.2] From $t=1$ to $|A|$, play each $a\in A$ once and obtain $X_{a,t}$ and $Y_{a,t}$ independently.
\item[1.3] Set $T_{a}(|A|) = 1$ for all $a\in A$ and $t=|A|+1$.
\end{itemize}
\item[2.] \textbf{Loop while $t\leq T$} 
\begin{itemize}
\item[2.1] Play $I_t^{\pi} \in \argmin_{a\in A} \left ( \min(\bar{\Delta}_a^\epsilon(t), \bar{\Phi}_a^\epsilon(t)) \sqrt{T_a(t)} \right )$.
\item[2.2] Obtain $X_{I^{\pi}_t,t}$ and $Y_{I^{\pi}_t,t}$ independently
and $T_{I^{\pi}_t}(t) \leftarrow T_{I^{\pi}_t}(t-1) + 1$ and $t\leftarrow t+1$.
\end{itemize}
\item[3.] \textbf{Output:}
\begin{itemize}
\item[3.1] Obtain $A_T^f(\epsilon) = \{a\in A | \bar{Y}_{T_a(T)} \leq C \}$ and $A_T^*(\epsilon) = \{a\in A | \bar{X}_{T_a(T)} \geq \mu^* \}$.
\item[3.2] Output $A_T^*(\epsilon) \cap A_T^f(\epsilon)$. 
\end{itemize}
\end{itemize}
\vspace{0.5cm}

\subsection{Asymptotic Optimality}
\label{subsec:opt}

To analyze the behavior of CAPT, we start with the definition of 
a set of approximately feasible arms: For a given $\kappa \in \Re$, $A^{\kappa}_f := \{a\in A| C_a \leq C + \kappa \}$. 
Given $\epsilon \geq 0$, any set $S$ in $\mathcal{P}(A)$ is referred to as an \emph{$\epsilon$-feasible set} of arms \emph{if} $A^{-\epsilon}_f \subseteq S \subseteq A^{\epsilon}_f$,
where $\mathcal{P}(A)$ is the power set of $A$.
An arm $a$ in $A$ is \emph{$\epsilon$-feasible} if $a$ is an $\epsilon$-feasible set.
We also define 
a set of competing (optimality-candidate) arms: For a given $\kappa \in \Re$, $A^{\kappa}_* := \{a\in A| \mu_a \geq \mu^* + \kappa \}$.
Given $\rho \geq 0$, any set $K$ in $\mathcal{P}(A)$ is referred to as a \emph{$\rho$-competing set} of arms \emph{if} $A^{\rho}_* \subseteq K \subseteq A^{-\rho}_*$.
An arm $a$ in a $\rho$-competing set is $\rho$-competing.
Note that a $\rho$-competing arm is \emph{not necessarily feasible}
and that for any given \emph{$\epsilon$-feasible set} $S$ and \emph{$\rho$-competing set} $K$,
$A^{-\epsilon}_f \cap A^{\rho}_* \subseteq S \cap K \subseteq A^{\epsilon}_f \cap A^{-\rho}_*$.
The set $S \cap K$ in the previous identity is said to be a \emph{$(\epsilon,\rho)$-optimal} set
and an arm $a$ in $S\cap K$ is \emph{$(\epsilon,\rho)$-optimal}.

In the sequel, we consider the case where $\epsilon=\rho$ 
and refer to an $(\epsilon,\epsilon)$-optimal set as just an $\epsilon$-optimal set.
An arm in an $\epsilon$-optimal set is $\epsilon$-optimal.
If $\epsilon=0$, the $0$-feasible set corresponds to $A_f$ and the $0$-competing set is equal to $\{a\in A| \mu_a \geq \mu^* \}$, and the intersection of the two sets is equal to the solution set of $\argmax_{a\in A_f}\mu_a$.

The theorem below states about a finite-time lower bound to the probability that $A_T^*(\epsilon)\cap A_T^f(\epsilon)$ produced by CAPT at $T$ in the \textbf{Output} step is an $\epsilon$-optimal set for some general conditions.
The bound is given in terms of a problem-complexity denoted by $H(\epsilon):=\sum_{a\in A}\min (\Delta_a^\epsilon,\Phi_a^\epsilon)^{-2}$. 
This complexity must be very intuitive:
The performance of CAPT depends on the sum of the degrees of the hardness of each action between the cost-feasibility problem and the reward-optimality problem.
Note that if $\epsilon=0$, $H(\epsilon)$ becomes infinity because $\Delta_a^\epsilon=0$ for some $a\in A^*_f$. 
If the problem contains an arm $a$ that satisfies the constraint by equality such that $\Phi_a^\epsilon=0$,
 $H(\epsilon)$ become infinity again.
Therefore, we exclude such cases by requiring that $\epsilon >0$ but can be arbitrarily close to zero.

The assumption that $T \geq 2|A|$ in the theorem statement is due to a technical reason:
Obviously, to make CAPT run, the condition that $T \geq |A|$ is necessary due to the \textbf{Initialization} step.
We further observe that
there always exists $a$ in $A$ such that $T_a(T)-1 \geq \frac{T-|A|}{H(\epsilon) \min(\Delta_a^\epsilon,\Phi_a^\epsilon)^2}$.
Suppose not. Then $T-|A| = \sum_{a\in A} (T_a(T) -1) < \sum_{a\in A} \frac{T-|A|}{H(\epsilon) \min(\Delta_a^\epsilon,\Phi_a^\epsilon)^2} = \frac{T-|A|}{H(\epsilon)} \sum_{a\in A} \frac{1}{\min(\Delta_a^\epsilon,\Phi_a^\epsilon)^2} = T-|A|$, which is a contradiction.
By $T \geq 2|A|$ then, we can fix an action $a$ that satisfies the bound of $T_a(T)-1 \geq \frac{T}{2H(\epsilon) \min(\Delta_a^\epsilon,\Phi_a^\epsilon)^2}$ and that has been played at least two times by $T$ and can use the inequality in ``cleaning up" some terms to eventually obtain a bound on $T_a(T)$.
In addition, we impose the condition that $X_{a,t}$ and $Y_{a,t}$ are in $[0,1]$ for any $a$ and $t$ for the better exposition.
\\
\begin{thm}
\label{thm1} Assume that the reward and the cost distributions associated with all arms in $A$ have the support in $[0,1]$.
Then for any $\epsilon > 0$ and $T\geq 2|A|$, the output $A_T^*(\epsilon) \cap A_T^f(\epsilon)$ by CAPT at $T$ satisfies
\[
\Pr\{ A^{\epsilon}_* \cap A^{-\epsilon}_f \subseteq  A_T^*(\epsilon)  \cap A_T^f(\epsilon)  \subseteq A^{-\epsilon}_* \cap A^{\epsilon}_f \} \geq 1 - 2|A|T e^{-T/16H(\epsilon)}.
\]
\end{thm}
%\\
Before presenting the proof, we remark that
the idea of the proof basically follows the reasoning in the proof of Theorem 2 by Locatelli \emph{et al.}~\cite{locatelli} since CAPT is built upon APT.
But the proof here requires the more thoughts due to the different index to be manipulated. We also polish some arguments given in~\cite{locatelli}.
In particular, the simpler Hoeffding inequality~\cite{hoeff} is applied in a place where
a lower bound to some probability is obtained instead of nonidentifiable ``Sub-Gaussian 
martingale inequality" referred by Locatelli \emph{et al.}
The lower bound with the term of $O(|A|T e^{-T})$ in our result is looser than the stated
lower bound with $O(\log T e^{-T})$ to a related probability by Locatelli \emph{et al.}
However, the arguments of Locatelli \emph{et al.}~for the tighter $\log$-bound seems 
\emph{incomplete} at the steps of applying the Union bound.
In fact, Wang and Ahmed~\cite{wang} provide a related result for the cost-feasibility problem 
that has the \emph{same order} of $O(|A|T e^{-T})$.
Their approach is within the context of the ``sample average approximation"~\cite{kleywegt}.
Thus the method is not index-based and not adaptive.
In our terms, they analyzed the probability that $\{a\in A| \bar{Y}_{T_a(T)} \leq C \}$ is 
an $\epsilon$-feasible set when each action in $A$ is played $N$ times equally, that is, 
$T_a(T) = N$ for all $a\in A$.
It is not clear how the adaptive index-based approach of APC makes a jump from
$O(T)$ to $O(\log T)$ in the order in Locatelli \emph{et al.}'s proof.
Besides, Locatelli \emph{et al.}'s theorem statement includes 
the case of $\epsilon=0$, which will lead to the non-asymptotic optimality.
%Note that the definition of the probability of ``making a mistake by $\epsilon$-tolerance" in~\cite{locatelli} is equivalent to the probability that $A_T^f(\epsilon)$ in CAPT is not a $\epsilon$-feasible set.
\\
\begin{proof}
Define an event $\xi$ such that with a given $\delta > 0$,
\begin{eqnarray*}
\lefteqn{\xi = \biggl \{ \forall a\in A, \forall T_a(T)\in \{1,2,...,T \},}\\
& & |\bar{X}_{T_a(T)} - \mu_a | \leq \sqrt{\frac{T\delta^2}{H(\epsilon) T_a(T)}}
\bigwedge |\bar{Y}_{T_a(T)} - C_a | \leq \sqrt{\frac{T\delta^2}{H(\epsilon) T_a(T)}}
\biggr \}.
\end{eqnarray*}

Fix any $a^*$ in $A$ such that $T_{a^*}(T)-1 \geq \frac{T}{2H(\epsilon) \min(\Delta_{a^*}^\epsilon,\Phi_{a^*}^\epsilon)^2}$ and
fix $t$ as the smallest $s$ in $\{|A|+1,...,T\}$ such that $T_{a^*}(s)=T_{a^*}(s+1),...,=T_{a^*}(T)$. In other words, $t$ is the last time $a^*$ was played and satisfies that $T_{a^*}(t) \geq T_{a^*}(T)-1$.

On $\xi$ we have that for all $i\in A$,
\[
     | \bar{X}_{T_i(t)} - \mu_i | \leq \sqrt{\frac{T\delta^2}{H(\epsilon) T_i(t)}}.
\] This implies that
for all $i\in A$,
\[
\Delta_i^\epsilon - \sqrt{\frac{T\delta^2}{H(\epsilon) T_i(t)}} \leq \bar{\Delta}_i^\epsilon(t) \leq \Delta_i^\epsilon + \sqrt{\frac{T\delta^2}{H(\epsilon) T_i(t)}}
\] 
because for all $i\in A$,
$
      | \bar{X}_{T_i(t)} - \mu_i | \geq |\bar{\Delta}_i^\epsilon(t) - \Delta_i^\epsilon |
$ 
where we recall $\Delta_i^\epsilon = |\mu_i - \mu^*| + \epsilon$ and $\bar{\Delta}_i^\epsilon(t) = |\bar{X}_{T_i(t)} - \mu^*| + \epsilon.$ 

Similarly, for all $i\in A$,
\[
\Phi_i^\epsilon - \sqrt{\frac{T\delta^2}{H(\epsilon) T_i(t)}} \leq \bar{\Phi}_i^\epsilon(t) \leq \Phi_i^\epsilon + \sqrt{\frac{T\delta^2}{H(\epsilon) T_i(t)}}
\] where we recall
$\Phi_i^\epsilon = |C_i - C| + \epsilon$ and $\bar{\Phi}_i^\epsilon(t) = |\bar{Y}_{T_i(t)} - C| + \epsilon.$

Because $a^*$ was played at $t$, $a^*$ achieves the value of the minimum index, i.e., $K_{a^*}(t) \leq K_i(t)$ for all $i\in A$. Recall that
\[
   K_{a^*}(t) = \min \left (\bar{\Delta}_{a^*}^\epsilon(t),\bar{\Phi}_{a^*}^\epsilon(t) \right ) \sqrt{T_{a^*}(t)}.
\]
From the two inequalities of 
$\Delta_{a^*}^\epsilon - \sqrt{\frac{T\delta^2}{H(\epsilon) T_{a^*}(t)}} \leq \bar{\Delta}_{a^*}^\epsilon(t)$ and
$\Phi_{a^*}^{\epsilon} - \sqrt{\frac{T\delta^2}{H(\epsilon) T_{a^*}(t)}} \leq \bar{\Phi}_{a^*}^\epsilon(t)$, 
it follows that 
\[
\min \left ( \Delta_{a^*}^\epsilon - \sqrt{\frac{T\delta^2}{H(\epsilon) T_{a^*}(t)}}, \Phi_{a^*}^\epsilon - \sqrt{\frac{T\delta^2}{H(\epsilon) T_{a^*}(t)}} \right )
\leq \min \left ( \bar{\Delta}_{a^*}^\epsilon(t), \bar{\Phi}_{a^*}^\epsilon(t) \right ).
\] Thus we have that
\[
\min( \Delta_{a^*}^\epsilon, \Phi_{a^*}^\epsilon ) - \sqrt{\frac{T\delta^2}{H(\epsilon) T_{a^*}(t)}} \leq \min \left ( \bar{\Delta}_{a^*}^\epsilon(t), \bar{\Phi}_{a^*}^\epsilon(t) \right ).
\] 
Multiplying both sides of the above inequality by $\sqrt{T_{a^*}(t)}$ and using $\sqrt{T_{a^*}(t)}
\geq \sqrt{T}/\sqrt{2H(\epsilon)\min(\Delta_{a^*}^\epsilon,\Phi_{a^*}^\epsilon)^2}$ and rearranging the terms 
leads to a lower bound to $K_{a^*}(t)$:
\begin{equation}
\label{eqn1}
      \left (\frac{1}{\sqrt{2}} - \delta \right ) \sqrt{\frac{T}{H(\epsilon)}}\leq K_{a^*}(t). 
\end{equation}
We now upper bound $K_i(t)$ for any $i\in A$. From the inequality for the bound of $\bar{\Delta}_i^\epsilon$, we have that
%\begin{eqnarray}
%\lefteqn{\hspace{-4cm} K_i(t) = \min \left ( \bar{\Delta}_i^\epsilon(t),\bar{\Phi_i}^\epsilon(t) \right ) \sqrt{T_i(t)} \leq \bar{\Delta}_i^\epsilon(t) \sqrt{T_i(t)}} \nonumber \\
%\hspace{4cm} & & \leq \left ( \Delta_i^\epsilon + \sqrt{\frac{T\delta^2}{H(\epsilon) T_i(t)}} \right ) \sqrt{T_i(t)}. \label{eqn2}
%\end{eqnarray}
%
\begin{equation}
K_i(t) = \min \left ( \bar{\Delta}_i^\epsilon(t),\bar{\Phi}_i^\epsilon(t) \right ) \sqrt{T_i(t)} \leq \bar{\Delta}_i^\epsilon(t) \sqrt{T_i(t)} \leq \left ( \Delta_i^\epsilon + \sqrt{\frac{T\delta^2}{H(\epsilon) T_i(t)}} \right ) \sqrt{T_i(t)}. \label{eqn2}
\end{equation}
Combining~(\ref{eqn1}) and~(\ref{eqn2}) results in
\[
    \left (\frac{1}{\sqrt{2}} - \delta \right ) \sqrt{\frac{T}{H(\epsilon)}} \leq   \Delta_i^\epsilon \sqrt{T_i(t)} + \delta \sqrt{\frac{T}{H(\epsilon)}}
\] for all $i\in A$. Rearranging the terms in the above inequality and from $T_i(T)\geq T_i(t)$,
\[
     (1-2\sqrt{2}\delta)^2 \frac{T}{2H(\epsilon)(\Delta_i^\epsilon)^2} \leq T_i(T).
\] In sum, $\xi$ implies that for any $i\in A$,
\begin{equation}
\label{eqn:capt}
   \mu_i - \Delta_i^\epsilon \times \frac{\sqrt{2}\delta}{1-2\sqrt{2}\delta} \leq \bar{X}_{T_i(T)} \leq
   \mu_i + \Delta_i^\epsilon \times \frac{\sqrt{2}\delta}{1-2\sqrt{2}\delta}.
\end{equation}
Set $\frac{\sqrt{2}\delta}{1-2\sqrt{2}\delta} = 1/2$ by letting $\delta=(4\sqrt{2})^{-1}$.
We show that the event $\xi$ implies that
\[
A^{\epsilon}_* \subseteq  \left \{j\in A| \bar{X}_{T_j(T))} \geq \mu^* \right \} \subseteq A^{-\epsilon}_*.
\]
%%reward-optimality case
For any $i\in A$ such that $\mu_i \geq \mu^* + \epsilon$, $\Delta_i^\epsilon = \mu_i-\mu^* + \epsilon$. By $\mu_i - \frac{1}{2} \Delta_i^\epsilon \leq \bar{X}_{T_i(T)}$,
\[
  \bar{X}_{T_i(T)} - \mu^* \geq \mu_i - \frac{1}{2}\Delta_i^\epsilon - \mu^* =  \mu_i - \frac{1}{2} (\mu_i-\mu^* + \epsilon) - \mu^*  \geq 0
\] making $\bar{X}_{T_i(T)} \geq \mu^*$. On the other hand, for any $i\in A$ such that
$\mu_i < \mu^* - \epsilon$, $\Delta_i^\epsilon = \mu^* - \mu_i + \epsilon$ and this results in
$\bar{X}_{T_i(T)} < \mu^*$.
%It follows that
%\[
%   A^{\epsilon}_* \subseteq  \left \{j| \bar{X}_{T_j(T))} \geq \mu^*, j\in A \right \} \subseteq A^{-\epsilon}_*.
%\]

%%%%feasibility case
We next consider the cost-feasibility case. By the same method as in~(\ref{eqn2}), on $\xi$ we have that for any $i\in A$,
\[
K_i(t) \leq \bar{\Phi}_i^\epsilon(t) \sqrt{T_i(t)} \leq \left ( \Phi_i^\epsilon + \sqrt{\frac{T\delta^2}{H(\epsilon) T_i(t)}} \right ) \sqrt{T_i(t)}.
\]
Then following the similar arguments as in the reward-optimality case leads to the inequality of
\[
   C_i - \Phi_i^\epsilon \times \frac{\sqrt{2}\delta}{1-2\sqrt{2}\delta} \leq \bar{Y}_{T_i(T)} \leq C_i + \Phi_i^\epsilon \times \frac{\sqrt{2}\delta}{1-2\sqrt{2}\delta}.
\] 
With $\delta=(4\sqrt{2})^{-1}$, 
for any $i\in A$ such that $C_i > C + \epsilon$, $\Phi_i^\epsilon = C_i  -C + \epsilon$ and
it follows that $\bar{Y}_{T_i(T)} > C$ because
\[
 \bar{Y}_{T_i(T))} - C \geq C_i - \frac{1}{2}\Phi_i^\epsilon - C =  C_i - \frac{1}{2} (C_i-C +\epsilon)  - C = \frac{1}{2}(C_i -C - \epsilon) > 0
\] by $C_i - \frac{1}{2}\Phi_i^\epsilon \leq \bar{Y}_{T_i(T)}$.
Furthermore, for any $i\in A$ such that $C_i \leq C - \epsilon$, $\Phi_i^\epsilon = C - C_i  + \epsilon$. 
Because
\[
\bar{Y}_{T_j(T)} - C \leq C_i + \frac{1}{2} \Phi_i^\epsilon - C =  C_i + \frac{1}{2}(C - C_i  + \epsilon)  - C 
= \frac{1}{2}(C_i - C + \epsilon) \leq 0,
\] $\bar{Y}_{T_j(T)} \leq C$.
It follows that
$
   A^{-\epsilon}_f \subseteq \left \{j\in A| \bar{Y}_{T_j(T))} \leq C\right \} \subseteq A^{\epsilon}_f.
$

Putting the reward-optimality and the cost-feasibility arguments together (by independence), $\xi$ implies that  
\[
A^{\epsilon}_* \subseteq  \{j\in A| \bar{X}_{T_j(T))} \geq \mu^* \} \subseteq A^{-\epsilon}_*
 \mbox{ and } 
A^{-\epsilon}_f \subseteq \{j\in A| \bar{Y}_{T_j(T))} \leq C \} \subseteq A^{\epsilon}_f
.
\]

By applying the Union bound (Boole's inequality) and Hoeffding inequality~\cite{hoeff}, the probability of $\xi$ is lower bounded as follows:
\begin{eqnarray*}
\lefteqn{\Pr(\xi) = 1-\Pr(\xi^c)}\\
& & \geq 1 - \sum_{a\in A}\sum_{T_a(T)=1}^T \Biggl ( \Pr \left \{ |\bar{X}_{T_a(T)} - \mu_a | > \sqrt{\frac{T\delta^2}{H(\epsilon) T_a(T)}} \right \} \\
& & \hspace{5.2cm} + \Pr \left \{ |\bar{Y}_{T_a(T)} - C_a | > \sqrt{\frac{T\delta^2}{H(\epsilon) T_a(T)}} \right \} \Biggr ) \\
& & \geq 1 - |A| T e^{-2T\delta^2/H(\epsilon)} - |A| T e^{-2T\delta^2/H(\epsilon)} =  1 - 2|A| T e^{-2T/16H(\epsilon)}.
\end{eqnarray*} 
\end{proof}

\section{CAPT with Estimation Algorithm}
\label{sec:relax}

In this section, we provide an algorithm \emph{in a general form} that replaces $\mu^*$ with $\mu^*(t)$ where $\mu^*(t)$ denotes the estimate of $\mu^*$ at $t$.
We call the algorithm ``CAPT with Estimation" (CAPT-E) and refer to it as $\pi'$ wherever possible.
We discuss two examples below for the estimation.

\subsection{Algorithm}

The procedure is the same as that of CAPT except that the role of $\mu^*$ is replaced by $\mu^*(t)$. In particular, $\bar{\Delta}_a^\epsilon(t)$ in CAPT is changed with $\bar{\Delta}_a^{\epsilon,*}(t)$ where $\bar{\Delta}_a^{\epsilon,*}(t) := |\bar{X}_{T_i(t)} - \mu^*(t) | + \epsilon$. The set $A_T^*(\epsilon)$ in the \textbf{Output} step of CAPT is also changed with the set $\{a\in A | \bar{X}_{T_a(T)} \geq \mu^*(t)\}$.
We abuse the notations used in the previous section.
\\
\noindent\textbf{The CAPT with Estimation (CAPT-E) algorithm} $\pi'$
\begin{itemize}
\item[1.] \textbf{Initialization of CAPT}
\item[2.] \textbf{Loop while $t\leq T$} 
\begin{itemize}
\item[2.2] Play $I_t^{\pi'} \in \argmin_{a\in A} \left ( \min(\bar{\Delta}_a^{\epsilon,*}(t), \bar{\Phi}_a^\epsilon(t)) \sqrt{T_a(t)} \right )$.
\item[2.3] Obtain $X_{I^{\pi'}_t,t}$ and $Y_{I^{\pi'}_t,t}$ independently and $T_{I^{\pi'}_t}(t) \leftarrow T_{I^{\pi'}_t}(t-1) + 1$ and $t\leftarrow t+1$.
\end{itemize}
\item[3.] \textbf{Output:}
\begin{itemize}
\item[3.1] Obtain $A_T^f(\epsilon) = \{a\in A | \bar{Y}_{T_a(T)} \leq C \}$ and $A_T^*(\epsilon) = \{a\in A | \bar{X}_{T_a(T)} \geq \mu^*(t) \}$.
\item[3.2] Output $A_T^*(\epsilon) \cap A_T^f(\epsilon)$. 
\end{itemize}
\end{itemize}
\vspace{0.5cm}

In the next section, we discuss a general sufficient condition that makes CAPT-E achieve the asymptotic optimality and some example methods for estimation.

\subsection{Convergence}
\label{subsec:conv}

We reassume that $\epsilon > 0$ and $T\geq 2|A|$. Fix $a^*$ in $\{a \in A| T_a(T) \geq \frac{T}{2H(\epsilon) \min(\Delta_a^\epsilon,\Phi_a^\epsilon)^2} \}$ and fix $t$ as the last time $a^*$ was played. Notice that $t\rightarrow \infty$ as $T\rightarrow \infty$.

Obviously, in order for CAPT-E to achieve the asymptotic optimality, the following condition is
sufficient:
the relative distance to the optimal value from the reward sample-mean of each action $a$ at the horizon $T$, $\Delta_a^{\epsilon,*}(T)$, approaches the true value $\Delta_{a}^{\epsilon}$ as $T$ approaches infinity.
More precisely,
if $\Delta_a^{\epsilon,*}(t) \rightarrow \Delta_{a}^{\epsilon}$ for all $a\in A$ as $T \rightarrow  \infty$, then the probability that the output
$A_T^*(\epsilon) \cap A_T^f(\epsilon)$ by CAPT-E at $T$ is an $\epsilon$-competing set converges to one as $T\rightarrow \infty$.

We argue now that the above statement is indeed true.
As in the proof of Theorem~\ref{thm1}, let us define an event $\xi$ (from CAPT-E) such that with 
given $\delta > 0$ and $\delta_f > 0$,
\begin{eqnarray*}
\lefteqn{\xi = \biggl \{ \forall a\in A, \forall T_a(T)\in \{1,2,...,T \},}\\
& & |\bar{X}_{T_a(T)} - \mu_a | \leq \sqrt{\frac{T\delta^2}{H(\epsilon) T_a(T)}}
\wedge |\bar{Y}_{T_a(T)} - C_a | \leq \sqrt{\frac{T \delta_f^2}{H(\epsilon) T_a(T)}}
\biggr \}.
\end{eqnarray*}
On $\xi$, because
$
     | \bar{X}_{T_i(t)} - \mu_i | \leq \sqrt{\frac{T\delta^2}{H(\epsilon) T_i(t)}}
$ for all $i\in A$,
$
      | \bar{X}_{T_i(t)} - \mu_i | \geq |\bar{\Delta}_i^{\epsilon,*}(t) - \Delta_i^{\epsilon,*}(t) |
$ for all $i\in A$
where we define $\Delta_i^{\epsilon,*}(t) := |\mu_i - \mu^*(t)| + \epsilon$ and $\bar{\Delta}_i^{\epsilon,*}(t) := |\bar{X}_{T_i(t)} - \mu^*(t)| + \epsilon.$ 
It follows that 
%\[
%\min \left ( \bar{\Delta}_i^{\epsilon,*}(t) - \sqrt{\frac{T\delta^2}{H(\epsilon) T_i(t)}}, \Phi_i^{\epsilon} - \sqrt{\frac{T\delta^2}{H(\epsilon) T_i(t)}} \right )
%\leq \min \left ( \bar{\Delta}_i^{\epsilon,*}(t), \bar{\Phi}_i^\epsilon(t) \right )
%\] 
%having that
\[
\min( \Delta_{a^*}^{\epsilon,*}(t), \Phi_{a^*}^\epsilon ) - \sqrt{\frac{T\delta^2}{H(\epsilon) T_{a^*}(t)}} \leq \min \left ( \bar{\Delta}_{a^*}^{\epsilon,*}(t), \bar{\Phi}_{a^*}^\epsilon(t) \right ).
\] Multiplying both sides by $\sqrt{T_{a^*}(t)}$ and using $\sqrt{T_{a^*}(t)} 
\leq \sqrt{T}/\sqrt{2H(\epsilon)\min(\Delta_{a^*}^\epsilon,\Phi_{a^*}^\epsilon)^2}$ and rearranging the terms lead to
\begin{equation}
      \left (\frac{1}{\sqrt{2}} \times \frac{\min \left ( \bar{\Delta}_{a^*}^{\epsilon,*}(t), \bar{\Phi}_{a^*}^\epsilon(t) \right )}{\min(\Delta_{a^*}^\epsilon,\Phi_{a^*}^\epsilon)} - \delta \right ) \sqrt{\frac{T}{H(\epsilon)}}\leq \kappa_{a^*}(t).
\end{equation} 

An upper bound on $\kappa_i(t)$ for $i\in A$ is obtained by
%\begin{eqnarray}
%\lefteqn{\hspace{-4cm} K_i(t) = \min \left ( \bar{\Delta}_i^{\epsilon,*}(t),\bar{\Phi_i}^\epsilon(t) \right ) \sqrt{T_i(t)} \leq \bar{\Delta}_i^{\epsilon,*}(t) \sqrt{T_i(t)}} \nonumber \\
%\hspace{4cm} & & \leq \left ( \Delta_i^{\epsilon,*}(t) + \sqrt{\frac{T\delta^2}{H(\epsilon) T_i(t)}} \right ) \sqrt{T_i(t)}.
%\end{eqnarray}
%
\[
\kappa_i(t) = \min \left ( \bar{\Delta}_i^{\epsilon,*}(t),\bar{\Phi}_i^\epsilon(t) \right ) \sqrt{T_i(t)} \leq \bar{\Delta}_i^{\epsilon,*}(t) \sqrt{T_i(t)} \leq \left ( \Delta_i^{\epsilon,*}(t) + \sqrt{\frac{T\delta^2}{H(\epsilon) T_i(t)}} \right ) \sqrt{T_i(t)}.
\]
Let 
\[f_{a^*}(t) = \frac{\min \left ( \bar{\Delta}_i^{\epsilon,*}(t), \bar{\Phi}_i^\epsilon(t) \right )}{\min(\Delta_{a^*}^\epsilon,\Phi_{a^*}^\epsilon)}.
\]
Combining the lower and the upper bounds, we have that for all $i\in A$,
\[
    \left (\frac{1}{\sqrt{2}} f_{a^*}(t) - \delta \right ) \sqrt{\frac{T}{H(\epsilon)}} \leq \Delta_i^{\epsilon,*}(t) \sqrt{T_i(t)} + \delta \sqrt{\frac{T}{H(\epsilon)}}
\] and rearranging the terms 
\[
     \sqrt{\frac{T}{H(\epsilon)}} 
     \times \frac{\left (\frac{f_{a^*}(t)}{\sqrt{2}} - 2\delta \right )}{\Delta_i^{\epsilon,*}(t)} \leq \sqrt{T_i(T)}.
\] Applying this bound on $T_i(T)$ to 
$
     | \bar{X}_{T_i(T)} - \mu_i | \leq \sqrt{\frac{T\delta^2}{H(\epsilon) T_i(T)}}
$ in $\xi$ leads to
\begin{equation}
\label{eqn:capt-e}
   \mu_i - \Delta_i^{\epsilon,*}(t) \times \frac{\sqrt{2}\delta}{f_{a^*}(t)-2\sqrt{2}\delta} \leq \bar{X}_{T_i(T)} \leq
   \mu_i + \Delta_i^{\epsilon,*}(t) \times \frac{\sqrt{2}\delta}{f_{a^*}(t)-2\sqrt{2}\delta}
\end{equation} for all $i\in A$.

At this point, we can see that as $T\rightarrow \infty$,
$(\ref{eqn:capt-e})$ approaches (\ref{eqn:capt}) that we derived for CAPT when $\mu^*$ is known.
Note that if $\Delta_i^{\epsilon,*}(t) \rightarrow \Delta_i^{\epsilon}$ for all $i\in A$, $f_{a^*}(t)$ also approaches one because $\epsilon > 0$.
More formally, for
every given $\eta_1, \eta_2>0$, there exists a finite $T(\eta_1,\eta_2) > t(\eta_1,\eta_2) > 0$ such that for all $T > T(\eta_1,\eta_2)$ and $t > t(\eta_1,\eta_2)$, $| \Delta_i^{\epsilon,*}(t) - \Delta_i^{\epsilon} | \leq \eta_1$ and $| f_{a^*}(t) - 1  | \leq \eta_2$.
With setting $\delta = (1 - \eta_2)/4\sqrt{2}$, for such $t >  t(\eta_1,\eta_2)$ and $T > T(\eta_1,\eta_2)$, $\xi$ implies that
\begin{equation}
\label{eqn:capt-e-err}
   \mu_i - \frac{1}{2}(\Delta_i^{\epsilon}+\eta_1) \leq \bar{X}_{T_i(T)} \leq \mu_i + \frac{1}{2}(\Delta_i^{\epsilon}+\eta_1).
\end{equation}
Because the cost-feasibility part is the same as CAPT's, with $\delta_f = (4\sqrt{2})^{-1}$, on $\xi$ we have that
\[
   A^{-\epsilon}_f \subseteq \left \{j\in A| \bar{Y}_{T_j(T)} \leq C \right \} \subseteq A^{\epsilon}_f.
\]

Because we can make $\eta_1$ (and $\eta_2$) arbitrarily close to zero and the probability of $\xi$ 
also converges to one for any $\delta, \delta_f > 0$, after a sufficiently large $T$,
$A_T^*(\epsilon) \cap A_T^f(\epsilon)$ reaches to the limit of an $\epsilon$-competing set.

\subsection{Example}
\label{subsec:exam}

The immediate question is then what approximation scheme makes the sufficient condition satisfiable that $\Delta_a^{\epsilon,*}(t) \rightarrow \Delta_a^{\epsilon}$ for all $a\in A$ as $T\rightarrow \infty$. (And if such a scheme is available, the next question would be about the convergence speed.)

The difficulties around estimating $\mu^*$ are that first, $\mu^*$ is the \emph{expected} value and second, it is not simply $\max_{a\in A}\mu_a$ but $\max_{a\in A_f} \mu_a$.
Estimating the optimal value ``efficiently" is indeed a challenging open problem not just in CMAB but also in unconstraint MAB.
It seems not easy to avoid the curse of the law of large numbers or the central limit theorem while estimating the optimal value. 

Still, the sample-average value approach similar to the sample average approximation~\cite{kleywegt}
would be the simplest and the most straightforward 
approach: We set
\[\mu^*(t) = \max_{ j\in \{ i \in A| \bar{Y}_{T_i(t)} \geq C \} } \bar{X}_{T_j(t)}
\] if $\{i\in A| \bar{Y}_{T_i(t)} \geq C\} \neq \emptyset$ and a (pre-determined or arbitrarily chosen) constant in [0,1], otherwise. The value of $\mu_i, i\in A$ is estimated by $\bar{X}_{T_i(t)}$.

As we discussed before, due to the relationship of 
\[
     \sqrt{\frac{T}{H(\epsilon)}} 
     \times \frac{\left (\frac{f_{a^*}(t)}{\sqrt{2}} - 2\delta \right )}{\Delta_i^{\epsilon,*}(t)} \leq \sqrt{T_i(T)},
\] we can make $T_i(t)$ eventually approach infinity for any $i$ as $T\rightarrow \infty$.
Similar to the result of the sample average approximation, 
by the law of large numbers then, $\mu^*(t)$ will converge to $\mu^*$ (in probability) where the error diminishes asymptotically, i.e., with an $O(1/\sqrt{T})$ rate after a sufficiently large $T$.

Another approach we can consider is adapting the definition of the expected regret for the unconstrained MAB model.
We can set, for example,
\[
\mu^*(t) = \sum_{a\in \{ i \in A| \bar{Y}_{T_i(t)} \geq C\} } \bar{X}_{T_a(t)} \left ( \frac{T_a(t)}{t} \right ).
\] 
Suppose that $A^*_f=\{a^*\}$. If CAPT-E is asymptotically optimal, $T_{a^*}(t)\rightarrow \infty$ almost surely (a.s) as $t\rightarrow \infty$. If CAPT-E can guarantee $T_{a^*}(t)/t \rightarrow 1$ a.s~(and $T_{a}(t)/t \rightarrow 0$ for all $a \in A \setminus A^*_f$), $\mu^*(t)$ will converge to $\mu^*$ in the limit a.s.
Various adaptations would be possible. But the convergence analysis or establishing a bound of $|\mu^*(t)-\mu^*|$ for the corresponding adaptation is beyond the scope of this note and left as a future topic.

\section{Concluding Remark}

The establishment of the existence of an asymptotically optimal index-based deterministic algorithm 
for CMAB problems and the performance result of the algorithm
is expected to be a notable theoretical step to understand the solvability and the complexity of CMAB. 
An efficient algorithm for estimating $\mu^*$ combined with CAPT-E would be a good
candidate algorithm for solving CMAB problems. Devising such an algorithm is a good future research work.
In addition,
investigating the theoretical results of CAPT and CAPT-E 
by some experimental studies and doing some performance-comparison studies with
other (heuristic) algorithms is an important future work.

A direct application of CMAB is for approximately solving CMDP problems when
a set of (heuristic) policies is available at some initial state. Each policy 
can be simulated over a sample path over a finite horizon starting from the 
initial state and this can be viewed as obtaining a
sample reward in CMAB by viewing each policy as an arm. Then a best feasible 
policy (with some approximation degree) would be found at the initial state.

\end{document}